\documentclass[12pt, onesided]{amsart}
\usepackage{tikz,tikz-cd, tkz-graph}

\usetikzlibrary{positioning}
\usepackage{pgf,pgfplots,pgfplotstable,pgfmath}
\pgfplotsset{compat=1.18}

\usepackage{calc} 
\usepackage{stmaryrd} 
\usepackage{bookmark}

\usepackage{amsthm,amssymb,amsmath}

\usepackage{enumerate, enumitem}
\usepackage{verbatim}

\usepackage{fullpage}

\newtheorem{thmB}{Theorem}

\newtheorem{mainthm}{Theorem}

\newtheorem*{mainres*}{Theorem}
\newtheorem{thm}{Theorem}[section]

\newtheorem{lemma}[thm]{Lemma}
\newtheorem{prop}[thm]{Proposition}
\newtheorem{corollary}[thm]{Corollary}
\newtheorem*{definition*}{Definition}
\newtheorem*{thm*}{Theorem}
\newtheorem*{prop*}{Proposition}
\newtheorem*{lemma*}{Lemma}

\newtheorem*{corol*}{Corollary}

\newtheorem{conj}{Conjecture}
\theoremstyle{remark} 
\newtheorem*{example*}{Example}
\newtheorem{example}[thm]{Example}

\newtheorem*{remark*}{Remark}

\newcommand{\spin}{\ifmmode{\rm Spin}\else{${\rm spin}$\ }\fi}
\newcommand{\spinc}{\ifmmode{{\rm Spin}^c}\else{${\rm spin}^c$\ }\fi}

\newcommand{\spincs}{\mathfrak s}

\title{On lens spaces which may embed in the complex projective plane}
\author{Giulia Carfora}
\date{September 2026}

\begin{document}

\begin{abstract}
Using an obstruction derived by Owens and relying on Donaldson's Diagonalisation Theorem, we obtain a complete classification of lens spaces of Lisca type $(2)$ or $(3)$ that are unobstructed from embedding in a homology complex projective plane. This classification is further refined for lens spaces $L(p,q)$ with $p$ even, using an obstruction due to Lidman--Moore--Vazquez based on work of Lin. This article is extracted from the author's PhD thesis.
\end{abstract}
\maketitle
\section{Introduction}
For $p > q>0$ coprime integers, define the lens space $L(p,q)$ as the result of $-p/q$ surgery on the unknot in $S^3$. In this article we present the following result about lens spaces in a homology $\mathbb{C}P^2$:
\begin{mainthm}\label{thm: main}
    Any lens space $L(m^2,q)$, up to diffeomorphism and orientation reversal, that embeds in a homology $\mathbb{C}P^2$ and is of Lisca type $(2)$ or $(3)$ must have $m$ and $q$ of one of the following forms:
\begin{enumerate}[label=\textnormal{(\arabic*)}]
    \item $m= (2s+3)k^2+1$ and $q= m^2 - (2s+3)(m-1)$, with $s \geq 0$ and $k \neq 0$ even; 
    \item $m=24k^2+12k+2$ and $q=(m+1)(8k^2+4k+1)$, with $k \neq 0$;   
    \item $m=81k^2 -72k +17$ and $q=9(m+1)$, with $k =1 $ or $k$ even; 
    \item $m=3k^2-1$ and $q=6k^4-3k^2$, with $k \neq 0$ even; 
    \item $m=18k +23$ and $q=(m+1)(4k+5)$, with $k \geq 0$; 
    \item $m=81k^2+36k+5$ and $q=9(m+1)$, with $k\neq 0$ even; 
    \item $m=(2t+5)k^2 -1$ and $q=(2m-1)k^2$, with $t \geq 0$ and $k \neq 0$ even; 
    \item $m =18k +17$ and $q= (2m-1)(2k+2)$, with $k \geq 0$;
    \item $m = 18k + 31$ and $q = (4k+7)(m-1)$, with $k \geq 1$; 
    \item $m=128k^2+48k+5$ and $q=4(32k^2+12k+1)^2$, with $k \neq 0$; 
    \item $m=128k^2-80k+13$ and $q=4(32k^2-20k+3)^2$ , with $k \in \mathbb{Z}$; 
    \item $m=18k +19$ and $q=162(k+1)$, with $k \geq 0$; 
    \item either of $L(7^2,31)$, $L(9^2,50)$, $L(31^2, 270)$ or $L(11^2,50)$.\newline
\end{enumerate}
\end{mainthm}
We now explain what is meant by Lisca type $(2)$ or $(3)$. 
The classification presented relies on the fact that a lens space embedding in $\mathbb{C}P^2$ splits off a rational homology $4$-ball. In particular, lens spaces which do not bound such a rational homology $4$-ball cannot embed in $\mathbb{C}P^2$. In \cite{Lisca2007}, Lisca proved that the slice-ribbon conjecture holds for $2$-bridge knots and obtained a characterisation of lens spaces bounding rational homology $4$-balls. We focus on the characterisation of lens spaces. 
\begin{thm}[Lisca \cite{Lisca2007}]\label{thm:Liscaclassification} Let $p >q>0$ be coprime natural numbers. Up to diffeomorphism, both orientation-preserving and orientation-reversing, a lens space $L(p,q)$ bounds a rational homology $4$-ball if and only if $p = m^2$ and $q$ equals one of the following:
\begin{enumerate}[label=\textnormal{(\arabic*)}]
    \item $mk - 1 $ with $m>k>0$ and $\gcd(m,k) \in \{1,2\}$;
    \item $d(m \pm 1)$, where $d > 1$ divides $2m \mp 1$;
    \item $d(m \pm 1)$, where $d > 1$ is odd and divides $m \pm 1$.\newline
\end{enumerate}    
\end{thm}
\begin{remark*}
        We refer to the lens spaces in Theorem \ref{thm:Liscaclassification} as being of Lisca type $(i)$, where $i \in \{1,2,3\}$.
\end{remark*}

It is known that any orientable 3-manifold embeds smoothly in 5-space \cite{HirschImbedding}. Naturally, one might ask if it is also true that any such 3-manifold embeds in 4-space as well: the answer to this question has been known to be negative for a long time. Indeed, a homological argument due to Hantzsche shows that lens spaces, the simplest manifolds in terms of handle decomposition after $S^3$, cannot embed in $S^4$ \cite{Hantzsche1938}. 
The study of smooth embeddings of lens spaces (and of orientable $3$-manifolds in general) into $4$-manifolds did not stop at $S^4$. For instance, Edmonds--Livingston in \cite{EdmondsLivingston} showed that every lens space embeds smoothly in $\#_{n}\mathbb{C}P^2$ for some $n \in \mathbb{N}$.
Achieving a lower bound on $n$ for different lens spaces is an interesting problem on its own. \newline
\indent
This article summarises part of the author's PhD thesis \cite{CarforaPhD}, and addresses obstructions to smooth embeddings of individual lens spaces in a homology $\mathbb{C}P^2$. The reason for focusing on $\mathbb{C}P^2$ is that it is a positive-definite manifold with a simple structure in terms of handle decomposition. Note that proving that a lens space is obstructed from embedding in $\mathbb{C}P^2$ implies that $n \geq 2$ in the aforementioned problem.
Following from Lisca's work \cite{Lisca2007}, one sees that the problem of embedding a single lens space in $\mathbb{C}P^2$ is equivalent to that of embedding a rational homology ball it bounds, thus leading us to an intersection with a wider issue in algebraic geometry. Let $(m_{1}, m_{2}, m_{3})$ be a \emph{Markov triple}, that is, a triple of positive integers forming a solution to $m_1^2 +m_2^2 + m_3^2 = 3 m_{1}m_{2}m_{3}$. Such a triple gives rise to an embedding 
\[\sqcup_{i=1}^{3}
 B_{m_{i},k_{i}} \hookrightarrow \mathbb{C}P^2, 
\]
where each $B_{m_{i},k_{i}}$ is a rational homology ball obtained as the Milnor fiber of a smoothing of the corresponding cyclic quotient singularity $\frac{1}{m_{i}^2}(1, m_{i}k_{i}-1)$  with $\gcd(m_{i},k_{i})=1$; $B_{m_{i},k_{i}}$ has boundary $L(m_{i}^2, m_{i}k_{i}-1)$. This embedding arises from the smoothing of three singular points in the weighted projective space $\mathbb{P}(m_{1}^2, m_{2}^2,m_{3}^2)$. Hacking--Prokhorov showed in \cite{HackingProkhorov2010} that any projective surface $X$ admitting a smoothing to $\mathbb{C}P^2$ is $\mathbb{Q}$-Gorenstein deformation equivalent to some $\mathbb{P}(m_{1}^2, m_{2}^2,m_{3}^2)$, thus giving rise to an embedding of disjoint rational homology balls as above. Evans--Smith in \cite{EvansSmith2018} generalised this result to the symplectic case. In \cite{Owens2020}, Owens provided smooth non-symplectic embeddings of such rational homology balls. This work was later extended to exhibit more disjoint triples of rational homology balls by Lisca--Parma in \cite{LiscaParma22}, and again by Golla--Owens in \cite{GollaOwensQHB}. The focus on triples of rational homology balls was motivated by a conjecture of Koll\'ar, presented in Section \ref{sec:future}.  \newline
\indent
The work mentioned up to now provides smooth embeddings in a homology $\mathbb{C}P^2$ of lens spaces of Lisca type $(1)$. This article addresses embeddings of lens spaces of Lisca type $(2)$ and $(3)$, using an obstruction derived by Owens in \cite{Owens2020}. This obstruction relies on Donaldson's Diagonalisation Theorem from \cite{DonaldsonDiagonalization1987}. A lens space that embeds in $\mathbb{C}P^2$ splits off a rational homology ball $B$, and a positive-definite manifold $M$ such that $H_{1}(M; \mathbb{Z})=0$ and $H_2(M; \mathbb{Z}) \cong \mathbb{Z}$. We construct a new positive-definite manifold $X$ by gluing $M$ and $P$, the canonical positive-definite plumbing bounded by $-L(p,q)$, along their common boundary.
From Donaldson's Diagonalisation Theorem, it follows that the intersection lattices of $M$ and $P$ embed into the standard lattice on $\mathbb{Z}^{n+1}$, where $n=b_{2}(P)$, and satisfy some additional requirements derived by Owens, presented in Proposition \ref{prop:obs}.
Finding a lattice embedding satisfying Proposition \ref{prop:obs} implies that the lens space $L(p,q)$ is unobstructed from embedding in $\mathbb{C}P^2$. In general, we must work with the intersection lattices of both the plumbing bounded by $L(p,q)$ and that bounded by $-L(p,q)$. If it is not possible to find an embedding as in Proposition \ref{prop:obs} for either lattice, then we conclude that $L(p,q)$ does not embed in $\mathbb{C}P^2$. This analysis yields the classification in Theorem \ref{thm:MainThm}.
The result is then refined in Theorem \ref{thm:dthm} for lens spaces of type $(2)$ or $(3)$ with $p$ even, using an obstruction by Lidman--Moore--Vazquez \cite{LidmanMooreVazquez} presented in Theorem \ref{thm:dinv}. The latter concerns $d$-invariants and restrictions of spin structures from appropriate $4$-manifolds. 
\begin{mainthm}\label{thm:dthm}
 Any lens space $L(m^2,q)$ that embeds in a homology $\mathbb{C}P^2$, and is of Lisca type $(2)$ or $(3)$ with $m$ even, must be, up to diffeomorphism and orientation reversal, either $L(26^2, 243)$, or $L((24k^2+12k+2)^2, 3(8k^2+4k+1)^2)$ with $k \neq 0$.
\end{mainthm}
Combining Theorems \ref{thm:MainThm} and \ref{thm:dthm} yields Theorem \ref{thm: main}. All the results appearing in this article only rely on the homological description of the spaces involved. Therefore, any discussion about embeddings in $\mathbb{C}P^2$ should be understood to mean in a homology $\mathbb{C}P^2$.\\
\textbf{Acknowledgements.} The author would like to thank Paolo Aceto, Philipp Bader, Jim Belk, Marco Golla, and Brendan Owens for helpful conversations and guidance.
\section{Set-up}\label{sec:set-up}
  It can be shown using homological arguments that if a lens space $L(p, q)$ smoothly embeds in $\mathbb{C}P^2$, it separates it into two sub-manifolds: a rational homology ball and a positive-definite 4-manifold, $M$, such that $H_{1}(M; \mathbb{Z}) = 0$ and $H_{2}(M;\mathbb{Z}) \cong \mathbb{Z}$ \cite{Owens2020}. Lens spaces that do not bound rational homology balls are immediately excluded from our consideration. 
  Lens spaces bounding rational homology balls are presented in Theorem \ref{thm:Liscaclassification}, and our analysis of which ones may embed in $\mathbb{C}P^2$ focuses on those of type $(2)$ and $(3)$. \newline
  \indent
Let us introduce a key notion. An \emph{integer lattice} is a free abelian group together with a bilinear pairing taking values in $\mathbb{Z}$. The pair given by the free abelian group $\mathbb{Z}^n$, generated by the standard orthonormal basis $e_{1}, \dots, e_{n}$, together with the standard scalar product $e_{i}\cdot e_{j} = \delta_{ij}$ (where $\delta_{ij}$ is the Kronecker delta) is known as the standard lattice on $\mathbb{Z}^n$.    
A lattice $\Lambda(a_{1}, \dots, a_{n})$ can also be defined starting from integers $a_{1}, \dots, a_{n}$. Take the free abelian group $\mathbb{Z}^n$ with generators $v_{1}, \dots, v_{n}$. Define the symmetric bilinear pairing as follows:
   
    \[
    v_{i} \cdot v_{j} := 
    \begin{cases}
        a_{i} & \text{if $i=j$;} \\
        -1  & \text{if $|i-j| = 1$;} \\
         0  & \text{if $|i-j| > 1$.} \\
    \end{cases}
    \]
    This lattice is also called a \emph{linear lattice} and the generators $v_{1}, \dots, v_{n}$ are referred to as its vertices. The reason behind this nomenclature is that the symmetric bilinear pairing can be recovered as the adjacency matrix of a weighted linear graph. This is a finite graph with vertices $v_{1}, \dots, v_{n}$ where each $v_{i}$, $2 \leq i \leq n-1$, is connected only to $v_{i-1}$ and $v_{i+1}$, while $v_{1}$ is only connected to $v_{2}$ and $v_{n}$ is only connected to $v_{n-1}$.\newline\indent
    Recall that any lens space $L(p,q)$ can be expressed as the boundary of a positive-definite linear plumbing with weights $a_1, \dots, a_{n}$ determined by the following Hirzebruch-Jung continued fraction expansion:
    \[
    \frac{p}{p-q} = [a_{1}, a_{2}, \dots, a_{n} ]^{-} :=   a_1 - \cfrac{1}{a_2 - \cfrac{1}{a_3-\cfrac{1}{\raisebox{1.2ex}{$\ddots$} \,\, - \cfrac{1}{a_n}}}}, 
    \] 
    where each $a_i \geq 2$. The intersection lattice of such a plumbing is the linear lattice $\Lambda_{p/(p-q)}= \Lambda(a_{1}, \dots, a_{n})$.
    
\section{A lattice obstruction}\label{sec:dobs}
We may assume that $\partial M = L(p,q)$ and $\partial B = -L(p,q)$. From the integer long exact sequence of the pair $(M, L(p,q))$, it can be deduced that the map $H_2(M) \to H_2(M, L(p,q)) \cong H^2(M) $ is multiplication by $\pm p$. This map also represents the pairing of the intersection lattice of $M$, which is positive definite, and we may denote it as $\Lambda_M = (p)$ by a slight abuse of notation.
Construct a new $4$-manifold $X = M \bigcup\limits_{L(p,q)}P$ called a rational ``blow-up'' of $\mathbb{C}P^2$. Here $P$ is the canonical positive-definite plumbing of spheres bounded by $-L(p,q) \cong L(p,p-q)$ and associated to the continued fraction expansion of $p/q$ with weights $a_{1}, \dots, a_{n} $, for some $n \in \mathbb{N}$. The manifold $X$ is a closed, oriented, smooth 4-manifold with positive-definite intersection lattice. Hence, by Donaldson's Diagonalisation Theorem, the intersection lattice of $X$ is equivalent to the standard lattice on $\displaystyle{\mathbb Z}^{n+1}$. This implies the existence of a lattice embedding 
\[
\Lambda_M \oplus \Lambda_P \hookrightarrow \mathbb Z^{n+1},
\]
on which Owens derived additional constraints in \cite[Proposition 3.2]{Owens2020}, readapted and presented below.
\begin{prop}[Owens \cite{Owens2020}]\label{prop:obs} 
Let $B$ be a rational homology ball bounded by a lens space $-L(p,q)$ and suppose it embeds smoothly in $\mathbb{C}P^2$. Let $X = M \bigcup\limits_{L(p,q)} P$ be the resulting positive-definite rational blow-up of $\mathbb{C}P^2$. Then there exists a finite-index intersection lattice embedding
\[
\Lambda_{M} \oplus \Lambda_{P} \hookrightarrow {\mathbb{Z}^{n+1}}
\]
such that each unit vector $e \in {\mathbb{Z}^{n+1}}$ has non-zero pairing with each of $\Lambda_{M}$ and $\Lambda_{P}$. Moreover, the image of the generator of $\Lambda_{M}$ is a primitive vector in $\mathbb{Z}^{n+1}$ of square $p$.    
\end{prop}
The proposition enables us to analyse all lens spaces of Lisca type $(2)$ and $(3)$ and obtain the following preliminary classification. Combining this with Theorem \ref{thm:dthm} yields Theorem \ref{thm: main}, which differs from the following in the parity of $k$ in some of the families.
\begin{thm}\label{thm:MainThm}
Any lens space $L(m^2,q)$, up to diffeomorphism, that embeds in a homology $\mathbb{C}P^2$ and is of Lisca type $(2)$ or $(3)$ must have $m$ and $q$ of one of the following forms:
\begin{enumerate}[label=\textnormal{(\arabic*)}]
    \item $m= (2s+3)k^2+1$ and $q= m^2 - (2s+3)(m-1)$, with $s \geq 0$ and $k \neq 0, \pm 1$; 
    \item $m=24k^2+12k+2$ and $q=(m+1)(8k^2+4k+1)$, with $k \neq 0$;   
    \item $m=81k^2 -72k +17$ and $q=9(m+1)$, with $k \in \mathbb{Z}$; 
    \item $m=3k^2-1$ and $q=6k^4-3k^2$, with $k \neq 0, \pm 1$; 
    \item $m=18k +23$ and $q=(m+1)(4k+5)$, with $k \geq 0$; 
    \item $m=81k^2+36k+5$ and $q=9(m+1)$, with $k\neq 0$; 
    \item $m=(2t+5)k^2 -1$ and $q=(2m-1)k^2$, with $t \geq 0$ and $k \neq 0, \pm 1$; 
    \item $m =18k +17$ and $q= (2m-1)(2k+2)$, with $k \geq 0$;
    \item $m = 18k + 31$ and $q = (4k+7)(m-1)$, with $k \geq 1$; 
    \item $m=128k^2+48k+5$ and $q=4(32k^2+12k+1)^2$, with $k \neq 0$; 
    \item $m=128k^2-80k+13$ and $q=4(32k^2-20k+3)^2$ , with $k \in \mathbb{Z}$; 
    \item $m=18k +19$ and $q=162(k+1)$, with $k \geq 0$; 
    \item either of $L(7^2,31)$, $L(9^2,50)$, $L(31^2, 270)$ or $L(11^2,50)$.    
\end{enumerate}
\end{thm}

The process leading to the proof of Theorem \ref{thm:Liscaclassification} in \cite{Lisca2007} also yields an explicit description of all continued fraction expansions, also known as ``strings'',  arising from the lens space families involved. The string classification presented by Lisca has been recast in \cite{AcetoGollaLarsonLecuona2020} by Aceto--Golla--Larson--Lecuona with some slight modifications, and we shall adopt it in this article as well. We use the notation $2^{[s]}$ (or $2^{[t]}$) to indicate $s$ (or $t$) consecutive vertices of weight $2$. The strings of interest are the following:   
    \begin{itemize}
        \item $B^{1}_{s,t} = [2^{[t]}, 3, s+2, t+2,3,2^{[s]}]^{-}$ with $s\geq t >0$, of type $(2)$;
        \item $B^{2}_{s,t} = [2^{[t]}, s+3, 2, t+2, 3, 2^{[s]}]^{-}$ with $s\geq 0$, $t >0$, of type $(3)$;
        \item $C^{1}_{s,t}=[t+2, s+2,3,2^{[t]}, 4, 2^{[s]}]^{-}$ with $s>0$, $t \geq 0$, of type $(3)$;
        \item $C^{2}_{s,t} = [t+2, 2,s+3, 2^{[t]},4,2^{[s]}]^{-}$ with $s, t \geq 0$, of type $(2)$;
        \item $C^{3}_{s,t} = [t+3, 2, s+3, 3, 2^{[t]},3,2^{[s]}]^{-}$ with $s,t \geq 0$, of type $(3)$;
        \item $D^{1}_{s,t} = [t+3, 3, 2^{[s]}, 3,2^{[t]},3, s+3]^{-}$ where $s \geq t \geq 0$, of type $(2)$;
        \item $D^{2}_{s,t} = [t+3, 2^{[s]}, 4, 2^{[t]},3,s+2]^{-}$ with $s$, $t \geq 0$, of type $(3)$.
    \end{itemize}
Given a string above, reading it in reverse corresponds to the continued fraction expansion $p/q^*$, where $q^*$ is the multiplicative inverse of $q\mod p$. As $L(p,q)$ and $L(p,q^*)$ are diffeomorphic lens spaces, it suffices to work with the strings as presented above. We may keep reverse strings in mind to understand the relations between the lens spaces considered. For each string we denote its reverse by adding an $r$: for example, $rB^{1}_{s,t}$ is the reverse of $B^{1}_{s,t}$ and so on. Observe that the relevant differences between the classifications in \cite{Lisca2007} and \cite{AcetoGollaLarsonLecuona2020} are due to the following relations: $B^{1}_{s,0} = rB^{2}_{0,s}$, $B^{1}_{s,t} = rB^{1}_{t,s}$, $C^{1}_{0,t} = C^{2}_{0,t}$.
Any lens space $L(p,q)$ is orientation-reversing diffeomorphic to $L(p,p-q)$, which is known as the \emph{dual} lens space. This notion of duality also translates to strings, and we say that a string $[\mathbf{b}]^{-} = p/q$ is dual or \emph{complementary} to $[\mathbf{c}]^{-} = p/(p-q)$ \cite{AcetoGollaLarsonLecuona2020}.
Dualities between the strings presented above are as follows: $B^{1}_{s+1,t+1}$ is dual to $D^{1}_{s,t}$; $B^{2}_{s,t+1}$ is dual to $D^{2}_{s,t}$; $C^{1}_{s,t}$ is dual to $rC^{3}_{t,s-1}$; $C^{2}_{s,t}$ is dual to $rC^{2}_{t,s}$.\newline\indent
We outline the strategy employed to prove Theorem \ref{thm:MainThm}. For each string $[a_{1}, \dots, a_{n}]^{-}$, we consider the related linear graph and map each vertex to a linear combination of orthogonal basis vectors for $\mathbb{Z}^{n+1}$ so that the pairing is preserved. This corresponds to finding a lattice embedding of $\Lambda_P$, where $P$ is the plumbing corresponding to the graph, into $\mathbb{Z}^{n+1}$.
By Lisca's work, a lattice embedding of $\Lambda_P$ into $\mathbb{Z}^{n}$, and thus $\mathbb{Z}^{n+1}$, exists but does not satisfy Proposition \ref{prop:obs}. 
To show that the lens space bounding $P$ is unobstructed from embedding in $\mathbb{C}P^2$ we must find a lattice embedding of $\Lambda_P$ into $\mathbb{Z}^{n+1}$ that satisfies Proposition \ref{prop:obs}; this translates to making sure that each basis vector of $\displaystyle\mathbb{Z}^{n+1}$ appears in at least one of the linear combinations that the vertices of the linear graph are mapped to, and that the orthogonal complement of $\Lambda_P$ in $\mathbb{Z}^{n+1}$ --- which corresponds to the embedding $\Lambda_M \hookrightarrow \mathbb{Z}^{n+1}$ --- is generated by a vector with no zero entries and with square equal to $p$, the numerator of $p/q = [a_{1}, \dots, a_{n}]^{-}$. If such an embedding is found, the lens space $-L(p,q) \cong L(p,p-q)$ associated to the string is unobstructed from embedding in $\mathbb{C}P^2$. Naturally, this conclusion holds for $L(p,q)$ as well. However, if for a given string it is not possible to construct such an embedding, no conclusions can be drawn immediately.
In order to deduce that the corresponding lens space cannot embed in $\mathbb{C}P^2$, we must work with its dual string as well and verify that it is not possible to construct the required lattice embedding either. All of the relevant strings and possible lattice embeddings were analysed up to lattice automorphism of $\mathbb{Z}^{n+1}$.
The argument for each string starts by providing a general argument under the condition that the parameters $s$ and $t$ are ``big enough'' (defined in each case); then we analyse the infinite families of strings where only one of the parameters is ``big enough''; finally we analyse the finite number of (sporadic) cases left using GAP \cite{GAP4} and its \verb+OrthogonalEmbeddings+ command.  The complete argument to prove Theorem \ref{thm:MainThm} can be found in \cite[Chapter 4]{CarforaPhD}.
\begin{example}[Theorem \ref{thm:MainThm}, Case $(1)$] The string $B^{2}_{s,t} = [2^{[t]}, s+3, 2, t+2, 3, 2^{[s]}]^{-}$, with $s\geq 0$ and $t >0$, corresponds to the continued fraction expansion of 
\[\frac{p}{q}= \frac{(2st+2s+3t+4)^{2}}{(2st+2s+3t+4)^{2}-(2s+3)(2st+2s+3t+3)},\]
and is dual to $D^{2}_{s,t}$. We seek to analyse lattice embeddings in $\mathbb{Z}^{s+t+5}$, whose orthonormal basis we denote by $e_{1}, \dots, e_{s+1}, f_{1}, \dots, f_{t+1}, a, b, c$. We begin by discussing the case when $s$, $t \geq 3$. Figure \ref{fig:B2_general} shows the only way (up to automorphism) to obtain a lattice embedding of the lattice generated by $B^2_{s,t}$ in $\mathbb{Z}^{s+t+5}$, in a way that respects the bilinear pairing. 
\begin{figure}[h]
    \centering 
\begin{tikzpicture}[node distance=1.7cm,
  sphere/.style={circle,draw,minimum size=2mm, fill = black, inner sep=0pt},
  every node/.style={font=\small}
]

\node[sphere, label=above:{$f_{1}-f_{2}$}] (u1) {};
\node[sphere,right=of u1, label=above:{$f_{t}-f_{t+1}$}] (u2) {};
\node[sphere,right=of u2, label=above:{$v$}] (u3) {};
\node[sphere,right=of u3, label=above:{$a+b$}] (u4) {};
\node[sphere,right=of u4, label=above:{$w$}] (u5) {};
\node[sphere,right=of u5, label=above:{$-e_{1} + a- b$}] (u6) {};
\node[sphere,right=of u6, label=above:{$e_{1}-e_{2}$}] (u7) {};
\node[sphere,right=of u7, label=above:{$e_{s}-e_{s+1}$}] (u8) {};

\draw[dashed] (u1) -- (u2);
\draw (u2)--(u3)--(u4)--(u5)--(u6)--(u7);
\draw[dashed] (u7) -- (u8);
\draw[decorate,decoration={brace,mirror,amplitude=6pt}]
    ([yshift=-15pt]u1.south) -- ([yshift=-15pt]u2.south)
    node[midway, yshift=-12pt]{\small $t$ vertices};

\draw[decorate,decoration={brace,mirror,amplitude=6pt}]
    ([yshift=-15pt]u7.south) -- ([yshift=-15pt]u8.south)
    node[midway, yshift=-12pt]{\small $s$ vertices};
    
\node[below=2pt of u1] {$2$};
\node[below=2pt of u2] {$2$};
\node[below=2pt of u3] {$s+3$};
\node[below=2pt of u4] {$2$};
\node[below=2pt of u5] {$t+2$};
\node[below=2pt of u6] {$3$};
\node[below=2pt of u7] {$2$};
\node[below=2pt of u8] {$2$};
\end{tikzpicture}
   \caption{Potential embedding of the lattice generated by $B^{2}_{s,t}$ in $\mathbb{Z}^{s+t+5}$ with $s\geq 3$, $t\geq 3$.}
    \label{fig:B2_general}
\end{figure}
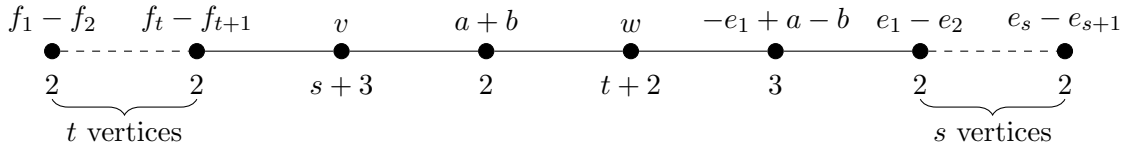
It remains to determine what the images $v$ and $w$ should be. Observe that the image of the vertex $v$ must be orthogonal to that of the vertices of weight 2 mapped to $e_i - e_{i+1}$; this implies that the integer coefficients of the $e_i$ vectors in the image of $v$ are equal. Similarly, the image of $v$ must be orthogonal to that of the vertices of weight 2 mapped to $f_i - f_{i+1}$, except for the vertex mapping to $f_t - f_{t+1}$, with which the pairing evaluates to $-1$. These types of arguments for $v$ and $w$ give rise to conditions on the linear combinations they are mapped to. We must have 
\begin{itemize}
    \item $v := \rho_{v} (e_{1} + \dots + e_{s+1}) + \eta_{v} (f_{1} + \dots + f_{t} ) + (\eta_{v} +1) f_{t+1} + \alpha_{v} a  -(1+\alpha_{v}) b + \gamma_{v} c$, and
    \item $w := \rho_{w} (e_{1} + \dots + e_{s+1}) + \eta_{w} (f_{1} + \dots + f_{t+1} ) + \alpha_{w} a -(1+\alpha_{w}) b + \gamma_{w} c$, 
\end{itemize}
where the additional parameters represent integer coefficients of the required linear combinations, subject to the conditions derived from the arguments mentioned above. We therefore need to solve the following system of equations: 
\begin{align*}
 \tag{1}  &v \cdot v =(s+1)\rho_{v}^2 +\eta_{v}^2 t + (\eta_{v} +1)^2 + \alpha_{v}^2 + (1+\alpha_{v})^2 + \gamma_{v}^2 = s+3 \implies |\rho_{v}| \leq 1   \\
 \tag{2} &w \cdot w =(s+1)\rho_{w}^2 + (t+1)\eta_{w}^2 + \alpha_{w}^2 + (1+\alpha_{w})^2 + \gamma_{w}^2 = t+2 \implies |\eta_{w}| \leq 1   \\
 \tag{3} & v \cdot (-e_1+a-b) =-\rho_{v} +2\alpha_{v} +1 =0 \implies \alpha_{v} = \frac{\rho_{v}-1}{2}\\
 \tag{4} &w \cdot (-e_1+a-b)= 2\alpha_{w} - \rho_{w} = -2 \implies \alpha_{w} =\frac{\rho_{w} - 2}{2} \\
 \tag{5} &v \cdot w = (s+1)\rho_{v}\rho_{w} +t\eta_{w}\eta_{v} + (\eta_{v} +1)\eta_{w} +\alpha_{v}\alpha_{w} + (1+\alpha_{v})(1+\alpha_{w}) +\gamma_{w}\gamma_{v} =0.
\end{align*}
Note that equation $(3)$ implies $\rho_{v} = \pm 1$, and equation $(4)$ implies $\rho_{w}$ is even.    
\begin{itemize}
    \item \textbf{Case $\rho_{v} = -1$:  } \newline
    After substituting for $\rho_{v}$ and simplifying, the system becomes:
    
    \begin{align*}
     \tag{1}  &\eta_{v}^2 t + (\eta_{v} +1)^2  + \gamma_{v}^2 = 1   \implies \eta_{v} = \gamma_{v} =0\\
     \tag{2} &(s+1)\rho_{w}^2 + (t+1)\eta_{w}^2 + \alpha_{w}^2 + (1+\alpha_{w})^2 + \gamma_{w}^2 = t+2 \implies |\eta_{w}| \leq 1   \\
     \tag{3} &\alpha_{v} = -1\\
     \tag{4} &\alpha_{w} =\frac{\rho_{w} - 2}{2} \\
     \tag{5} &(s+1)\rho_{w} + \alpha_{w} = \eta_{w} \implies \rho_{w} = \frac{2(\eta_{w}+1)}{2s+3}.
    \end{align*}

We are working under the assumption that $t \geq 3$. This implies that in equation $(1)$ $\eta_{v} =0$, and thus $\gamma_{v} =0$. A direct calculation shows that whenever $\eta_{w}$ equals 0 or 1 the system is inconsistent. This leaves us with the case $\eta_{w} = -1$. It immediately follows that $v = - (e_{1} + \dots + e_{s+1}) + f_{t+1} - a$ and $w =  -(f_{1} + \dots + f_{t+1} ) - a$. The resulting lattice embedding is Lisca's from \cite{Lisca2007}, which does not satisfy Proposition \ref{prop:obs}.

    \item \textbf{Case $\rho_{v} = 1$:  } \newline
    After substituting for $\rho_{v}$ and simplifying, the system becomes:
    
    \begin{align*}
     \tag{1}  &\eta_{v}^2 t + (\eta_{v} +1)^2  + \gamma_{v}^2 = 1   \implies \eta_{v} = \gamma_{v} =0\\
     \tag{2} &(s+1)\rho_{w}^2 + (t+1)\eta_{w}^2 + \alpha_{w}^2 + (1+\alpha_{w})^2 + \gamma_{w}^2 = t+2 \implies |\eta_{w}| \leq 1   \\
     \tag{3} &\alpha_{v} = 0\\
     \tag{4} &\alpha_{w} =\frac{\rho_{w} - 2}{2} \\
     \tag{5} &(s+1)\rho_{w} + \alpha_{w} + \eta_{w} +1 =0 \implies \rho_{w} = -\frac{2\eta_{w}}{2s+3}.
    \end{align*}

    As in the previous case, $t \geq 3$ implies that both $\eta_{v}$ and $\gamma_{v}$ are equal to 0. Note that $\eta_{w} = \pm 1$ leads to an inconsistent system. The case $\eta_{w} = 0$ directly leads to the mapping $v = e_{1} + \dots + e_{s+1} + f_{t+1}  - b $ and $w = -a +\gamma_{w} c$, whenever $t = \gamma_{w} ^2 -1$. This gives rise to a lattice embedding that employs each basis vector for $\mathbb{Z}^{s+t+5}$. The orthogonal complement is generated by the vector
    \[
        (  \underbrace{\makebox[2.4cm]{$2\gamma_{w},\dots, 2\gamma_{w}$}}_{\text{$s+1$ times}},
        \underbrace{\makebox[5.3cm]{$-\gamma_{w}(2s+3),\dots, -\gamma_{w}(2s+3)$}}_{\text{$t+1$ times}},
        \gamma_{w}, -\gamma_{w}, 1).
    \]
Note that $t \geq 3$ and $t = \gamma_{w}^2 -1$ imply that $\gamma_{w} \leq -2$ or $\gamma_{w} \geq 2$. In particular, this means that the vector presented has no zero entries. The square of this vector is $((2s+3)\gamma_{w}^2+1)^2 $, which is equal to $p$ when $t = \gamma_{w}^2 -1$. In this case, note that $q = ((2s+3)\gamma_{w}^2+1)^2 - (2s+3)^2\gamma_{w}^2$.
Thus, by Proposition \ref{prop:obs}, lens spaces with $p/q= B^2_{s, t}$ where $t = \gamma_{w}^2 -1 \geq 3$ are unobstructed from embedding in $\mathbb{C}P^2$.
\end{itemize}
To complete the proof of Case $(1)$ of Theorem \ref{thm:MainThm}, one needs to repeat the same type of argument for each $B^2_{s,t}$ string when either $s$ or $t$ is strictly smaller than 3.
\end{example}
\begin{example}[Theorem \ref{thm:MainThm}, Case $(13)$, $L(7^2,31)$] \label{ex: sporadic}
 Consider $p/q=B^2_{0,1} = [2, 3, 2,3,3 ]^{-}$ and the lattice embedding of $\Lambda_{p/q}$ in $\mathbb{Z}^6$. Input the matrix describing the pairing on $\Lambda_{p/q}$ into the \verb+OrthogonalEmbeddings+ command, together with the dimension of $\mathbb{Z}^n$ we want to embed in:
\begin{verbatim}
    OrthogonalEmbeddings([[2,-1,0,0,0],
                          [-1,3,-1,0,0],
                          [0,-1,2,-1,0],
                          [0,0,-1,3,-1],
                          [0,0,0,-1,3]], 6);
\end{verbatim}
This will output:
\begin{verbatim}
    
rec( norms := [ 1, 1, 1, 1, 1, 48/49, 48/49, 47/49, 47/49, 45/49, 41/49,
40/49, 40/49, 38/49, 34/49, 34/49, 31/49, 31/49, 27/49, 26/49, 24/49, 19/49 ],
    solutions := [ [ 1, 2, 3, 4, 5 ], [ 6, 7, 10, 12, 13, 21 ] ], 
    vectors := [ [ 1, -1, 0, 1, 0 ], [ -1, 0, 0, 1, 0 ], [ 0, -1, 1, -1, 1 ],
            [ 0, 1, 0, 0, 1 ], [ 0, 0, -1, 0, 1 ], [ 1, -1, 0, 0, 1 ], 
            [ -1, 0, 0, 0, 1 ], [ 0, 1, -1, 0, 1 ], [ 0, 0, 1, -1, 1 ],
            [ 0, 1, 0, -1, 1 ], [ 0, -1, 0, 0, 1 ], [ 0, 1, -1, 1, 0 ],
            [ 0, 0, 1, 0, 0 ], [ 0, -1, 0, 1, 0 ], [ 0, -1, 1, 0, 0 ], 
            [ 0, 0, -1, 1, 0 ], [ -1, 1, 0, 0, 0 ], [ 1, 0, 0, 0, 0 ],
            [ 0, 0, 0, -1, 1 ], [ 0, 1, 0, 0, 0 ], [ 0, 0, 0, 1, 0 ],
            [ 0, 0, 0, 0, 1 ] ] ).

\end{verbatim}
The list \verb+solutions+ contains two sublists, corresponding to two possible lattice embeddings in $\mathbb{Z}^6$. The elements of each sublist specify which elements should be taken from \verb+vectors+ as the rows of a matrix describing the embedding. Note that the first embedding in \verb+solutions+ employs five vectors: this is Lisca's embedding found in \cite{Lisca2007}, not satisfying Proposition \ref{prop:obs}. The second list contains six vectors and specifies the embedding in Figure \ref{fig:B2_0_1}.
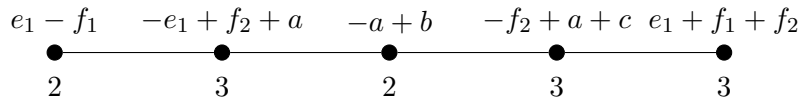
\begin{figure}[h]
    \centering 
\begin{tikzpicture}[node distance=2cm,
  sphere/.style={circle,draw,minimum size=2mm, fill = black, inner sep=0pt},
  every node/.style={font=\small}
]

\node[sphere, label=above:{$e_{1} - f_{1}$}] (u1) {};
\node[sphere,right=of u1, label=above:{$-e_{1} + f_{2}+a$}] (u2) {};
\node[sphere,right=of u2, label=above:{$-a+b$}] (u3) {};
\node[sphere,right=of u3, label=above:{$-f_{2} +a+c$}] (u4) {};
\node[sphere,right=of u4, label=above:{$e_{1}+f_{1}+f_{2}$}] (u5) {};

\draw (u1) -- (u2)--(u3)--(u4)--(u5);

\node[below=2pt of u1] {$2$};
\node[below=2pt of u2] {$3$};
\node[below=2pt of u3] {$2$};
\node[below=2pt of u4] {$3$};
\node[below=2pt of u5] {$3$};

\end{tikzpicture}
   \caption{Embedding of the lattice generated by $B^{2}_{0,1}$ in $\mathbb{Z}^{6}$.}
    \label{fig:B2_0_1}
\end{figure}
The orthogonal complement is generated by the vector 
\[
(1,1,-2,3,3,-5),
\]
which has square equal to $49 = p$. Thus, the lens space $L(49, 31)$ is unobstructed from embedding into  $\mathbb{C}P^2$.
\end{example}

\section{A $d$-invariant obstruction}\label{sec: lmvobs}
In this section, we refine the result presented in Theorem \ref{thm:MainThm} via an obstruction due to Lidman--Moore--Vazquez \cite{LidmanMooreVazquez}, after work by Lin \cite{LinPin2}. 
In the context of developing Floer homology in \cite{OzsvathSzabo2003AbsolutelyGraded}, Ozsv\'{a}th--Szab\'{o} introduced a ``correction term'', called the $d$-invariant $d(Y, \spincs)$, which is the rational number corresponding the minimal degree of any non-torsion class in $HF^{+}(Y, \spincs)$, associated to a rational homology sphere $Y$ with a $\operatorname{spin}^{\mathbb{C}}$ structure $\spincs$. 
An important property of $d$-invariants for us is that given $(Y, \spincs)$ as above, we have that
\[
d(-Y, \spincs) = -d(Y, \spincs),\]
where $-Y$ has opposite orientation to $Y$ \cite{OzsvathSzabo2003AbsolutelyGraded}.\newline\indent
Let us introduce the key concepts for our discussion. A compact, orientable, smooth $4$-manifold $W$ is a \emph{Koll\'ar manifold} if $H_{1}(W;\mathbb{Z})=0$, $H_2(W;\mathbb{Z}) \cong \mathbb{Z}$, and $\partial W$ is a union of rational homology $3$-spheres \cite{GollaOwensTopBMY}. 

\begin{prop}[Golla--Owens \cite{GollaOwensTopBMY}]\label{thm:evenp}
    Let $B$ be a rational homology 4-ball embedded in $X$, a homology $\mathbb{C}P^2$, such that $H_1(\partial B; \mathbb{Z})$ is cyclic of even order. Then $M = X \backslash \mathring{B}$ is a spin Koll\'ar manifold. Moreover, the restriction of the unique spin structure on $M$ to $\partial M$ does not extend to $B$.
\end{prop}
Note that not only does a Koll\'ar manifold as above admit a spin structure, but it admits a unique one. The next result concerns $L$-spaces, of which lens spaces are an example.

\begin{thm}[Lidman--Moore-Vazquez \cite{LidmanMooreVazquez}]\label{thm:dinv} Let $M$ be a positive-definite, spin Koll\'ar manifold, and let $Y = \partial M$, with induced spin structure $\spincs$. If $Y$ is an $L$-space then $d(Y,\spincs)=-\frac{1}{4}$.
\end{thm}

The previous results can be combined as follows.

\begin{figure}
    \centering

    \begin{minipage}{0.45\textwidth}
          \centering
           \begin{tikzpicture}
          \node[anchor = south west, inner sep=0] at (0,0)
        {\includegraphics[scale=1.5]{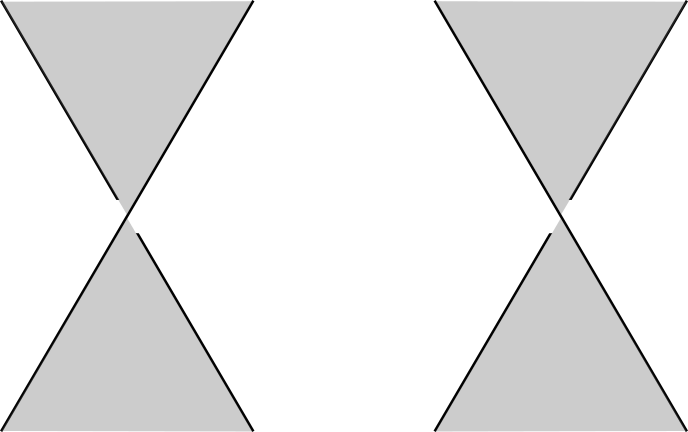}}; 
        \node at (1,-0.7) {$\eta(D) = +1$};
        \node at (4.4,-0.7) {$\eta(D) = -1$};
        \label{fig:crossings}
         \end{tikzpicture}
        \caption{Signed crossings.}
    \end{minipage}\hfill
    \begin{minipage}{0.45\textwidth}
        \centering
        \begin{tikzpicture}
        \node[anchor = south west, inner sep=0] at (0,0)
        {\includegraphics[scale=1.5]{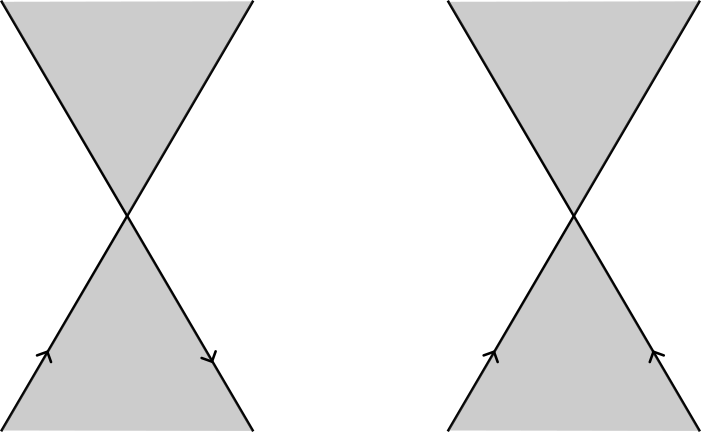}}; 
         \node at (1,-0.7) {Type I};
        \node at (4.4,-0.7) {Type II};
        \label{fig:typecrossings}
        \end{tikzpicture}
         \caption{Type I and II crossings.}
    \end{minipage}


\end{figure}

\begin{corollary}\label{cor:lmv}
    Let $m$ be even, and let $B$ be a rational homology 4-ball bounded by $-L(m^2,q)$ and embedded in $X$, a homology $\mathbb{C}P^2$. Then $M = X \backslash \mathring{B}$ is a positive-definite spin Koll\'ar manifold with induced spin structure $\spincs$ on the boundary $L(m^2,q)$, and $d(L(m^2, q),\spincs)=-\frac{1}{4}$.
\end{corollary}
A lens space $L(p,q)$ embedding in $\mathbb{C}P^2$ splits off a positive-definite Koll\'ar manifold $M$. Restricting ourselves to lens spaces from Theorem \ref{thm:MainThm} with $p=m^2$ even, we can make use of the obstruction from Corollary \ref{cor:lmv}. 
If we find a spin structure such that $d(L(p,q), \spincs) = -\frac{1}{4}$, then the conditions of Corollary \ref{cor:lmv} are satisfied, and we conclude that  $L(p,q)$ is unobstructed from embedding in $\mathbb{C}P^2$. Failure to satisfy the conditions of Corollary \ref{cor:lmv} for every spin structure implies that $L(p,q)$ cannot embed in $\mathbb{C}P^2$.
\newline\indent
Let us now discuss how to carry out the required computation.
A lens space $L(p,q)$ can be described as the double branched cover $\Sigma_2(S^3,L)$ of $S^3$ branched over a 2-bridge knot or link $L$. As we require $p$ to be even, here $L$ will always indicate a link. To compute $d(L(p,q), \spincs)$ we rely on the bijection established by Turaev between the quasi-orientations of a link $L \subset S^3$ and the spin structures on $\Sigma_2(S^3,L)$ \cite{Turaev}. A quasi-orientation of a link is an orientation up to overall reversal; thus, any link with $n$ components admits $2^{n-1}$ quasi-orientations.
We make use of the following result. 
\begin{lemma}[Donald--Owens \cite{DonaldOwens}] Let $L$ be a non-split alternating link with a quasi-orientation $o$, and let $\spincs_{o}$ be the corresponding spin structure on $\Sigma_2(S^3,L)$. Then 
    \[
d(\Sigma_2(S^3,L), \spincs_{o}) = -\frac{\sigma(L, o)}{4}.
\]
\end{lemma}
Here, $\sigma(L, o)$ denotes the signature of the link $L$. According to Gordon--Litherland \cite{Linksignature}, it can be computed by assigning a chessboard colouring to the diagram of an oriented link, that is, by colouring each region either black or white, so that any two adjacent regions have different colours. Let $R_0, \dots, R_n$ denote the white regions, and assign an incidence number $\eta(D) = \pm 1$ to each double point $D$ in the diagram as in Figure \ref{fig:crossings}. For $0 \leq i$ and $j \leq n$ define
\[
g_{ij} :=
\begin{cases}
    -\sum\eta(D) &  \text{summed over double points incident to $R_i$, $R_j$, with $i \neq j$}\\
    -\sum\limits_{k=0, \dots, n,\,  k \neq i} g_{ik} & \text{if $i=j$}.\\
\end{cases}
\]
The \emph{Goeritz matrix} associated to $L$ is the $n \times n$ symmetric matrix $G(L)= (g_{ij})$, with $i, j =1, \dots , n$. Then the signature of the link $L$, with quasi-orientation $o$, is given by $\sigma(L, o)= \sigma(G(L)) - \mu$, where $\mu = \sum \eta (D)$ summed over all double points of type II, shown in Figure 5.2. 
We may now compute $d(L(p,q), \spincs_{o})$ for any lens space with $p=m^2$ even and any choice of quasi-orientation $o$. As we are working with a 2-bridge link of two components, we need to consider two possible quasi-orientations, corresponding to two different spin structures $\spincs_{o}$ on $L(p,q)$. The computation was carried out for all lens spaces from the table of Theorem \ref{thm:MainThm} which admit even $p$. It produced the following result.

\begin{thmB}
 Any lens space $L(m^2,q)$ that embeds in a homology $\mathbb{C}P^2$, and is of Lisca type $(2)$ or $(3)$ with $m$ even, must be, up to diffeomorphism and orientation reversal, either $L(26^2, 243)$, or $L((24k^2+12k+2)^2, 3(8k^2+4k+1)^2)$ with $k \neq 0$.
\end{thmB}
It is interesting to observe that the lens spaces found in the theorem are all of type $(3)$ with $q=d(m+1)$, with string type $p/q = C^1_{s,t}$. Let us discuss the strategy employed and some subtleties due to the orientation of the lens spaces considered. Focus on the lens spaces with $p$ even from the list of Theorem \ref{thm:MainThm}. Recall that for a lens space $Y=L(p,q)$, we have assumed up to now that $\partial M= Y$ and $\partial B = -Y$. When computing the $d$-invariant, we will be interested in both the cases when $d(Y, \spincs_o) = \pm \frac{1}{4}$. This is to account for the assumption $\partial M= -Y$ and $\partial B = Y$ as well.
Therefore, we combine the approach employed to prove Theorem \ref{thm:MainThm} with the current one.
In particular, we need to pay special attention to some cases. 
To begin with, suppose that for $L(p,q)$ with corresponding 2-bridge link $L$ and quasi-orientations $o_1$ and $o_2$, neither $d(L(p,q), \spincs_{o_i})$ are $\pm \frac{1}{4}$; we conclude that $L(p,q)$ is obstructed from embedding in $\mathbb{C}P^2$ by Corollary \ref{cor:lmv}.
For the second case, suppose that $d(L(p,q), \spincs_{o_1}) = \frac{1}{4}$ and $d(L(p,q), \spincs_{o_2}) \neq -\frac{1}{4}$, or $d(L(p,q), \spincs_{o_1})=d(L(p,q), \spincs_{o_2}) =  \frac{1}{4}$: in this case we want to consider the lattice $\Lambda_{p/(p-q)}$ and verify whether it is possible to find an embedding satisfying Proposition \ref{prop:obs}. 
The third case to consider is when $d(L(p,q), \spincs_{o_1}) = -\frac{1}{4}$ and $d(L(p,q), \spincs_{o_2}) \neq \frac{1}{4}$, or $d(L(p,q), \spincs_{o_1}) =d(L(p,q), \spincs_{o_2}) = - \frac{1}{4}$: in this case we want to consider the lattice $\Lambda_{p/q}$ and verify whether it is possible to find an embedding satisfying Proposition \ref{prop:obs}. In either case, failure to find such an embedding implies that $L(p,q)$ does not embed in $\mathbb{C}P^2$. 
The complete proof of Theorem \ref{thm:dthm} can be found in \cite[Chapter 5]{CarforaPhD}.
\begin{example}
 Consider $L(p,q)$ with $p/q= C^1_{4k^2+2k-1,1} = [3, 4k^2+2k+1, 3, 2, 4, 2^{[4k^2+2k-1]}] $, with $k \neq 0$.
The relevant 2-bridge link is shown in Figure \ref{fig:C1_s_1_lidman}. The total number of crossings is $8k^2+4k+8$, and the signature of the associated Goeritz matrix $G$ is $\sigma(G(L))= 4k^2+2k+4$.
\begin{figure}[h]
    \centering
    \begin{tikzpicture}
    \node[anchor = south west, inner sep=0] at (0,0)
    {\includegraphics[scale = 0.8]{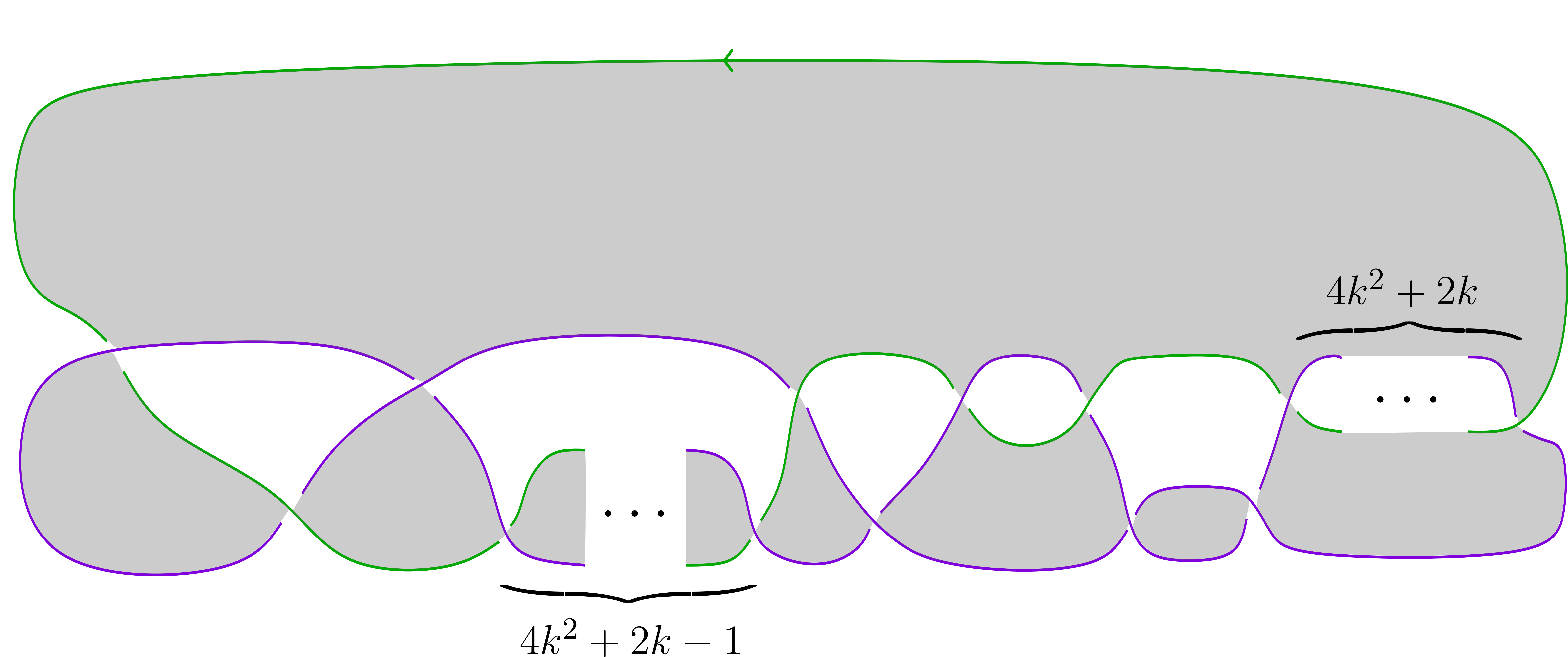}};

    \end{tikzpicture}
    
    \caption{The 2-bridge link associated with $\frac{p}{p-q}$ such that $\frac{p}{q}= C^1_{4k^2+2k-1,1}$.}
    \label{fig:C1_s_1_lidman}
\end{figure}
Let us compute the $d$-invariants for both choices of quasi-orientations of the link in Fig.~\ref{fig:C1_s_1_lidman}. One choice, $o_1$, yields $4k^2+2k+1$ type II crossings, and hence 
\[
d(L(p,p-q), \spincs_{o_1}) =-\frac{3}{4}.
\]
The second choice of quasi-orientation, $o_2$, yields $4k^2+2k+5$ type II crossings, and thus
\[
d(L(p,p-q), \spincs_{o_2}) = \frac{1}{4}.
\]
In light of our earlier discussion, we need to consider the lattice on $p/q=C^1_{4k^2+2k-1,1}$ and verify whether it is possible to find an embedding satisfying Proposition \ref{prop:obs}. The required embedding was found in \cite[Chapter 5]{CarforaPhD}, hence lens spaces $L(p,q)$ with $p=(24k^2+12k+2)^2$, $q=3(8k^2+4k+1)^2$, and $k \neq 0$ are unobstructed from embedding in $\mathbb{C}P^2$.
\end{example}

\section{Future directions}\label{sec:future}
In \cite{LiscaSum}, Lisca extended the work from \cite{Lisca2007} and classified which connected sums of lens spaces bound rational homology balls. Let $\mathcal{F}_{n}$ be the set
\[
\mathcal{F}_n = \left\{ \frac{m^2n}{mnk+1} \, \lvert \, m>k>0,\, \gcd(m,k)=1 \right\} \subset \mathbb{Q}, \quad n \geq 2. 
\]

\begin{thm}[Lisca \cite{LiscaSum}] Let $Y$ be an oriented 3-manifold homeomorphic to a connected sum of lens spaces. Then $Y$ smoothly bounds a rational homology ball if and  only if $Y$ is orientation-preserving homeomorphic to a connected sum in which each summand is (possibly orientation-reversing) homeomorphic to one of the following oriented $3$-manifolds:
\begin{enumerate}[label=\textnormal{(\arabic*)}]
    \item $L(p,q)$ with $p$ and $q$ as in Theorem \ref{thm:Liscaclassification};
    \item $L(p,q)\#L(p,p-q)$, with $\frac{p}{q} > 1$;
    \item $L(p_1,q_1) \# L(p_2,q_2)$, with $\frac{p_i}{q_i} \in \mathcal{F}_2$, for $i=1,2$;
    \item $L(p,q)\#L(n,n-1)$, with $\frac{p}{q} \in \mathcal{F}_{n}$ for some $n \geq 2$;
    \item $L(p_1,q_1)\# L(p_2, p_2-q_2)$, with $\frac{p_i}{q_i} \in \mathcal{F}_{n}$, $i=1,2$, for some $n \geq 2$.
\end{enumerate}
\end{thm}

A continuation of the work presented up to now would be to apply Donaldson's Diagonalisation Theorem to classify when a manifold $Y$ as above is unobstructed from embedding in $\mathbb{C}P^2$. Indeed, a connected sum of lens spaces $\#_{i}L(p_{i}, q_{i})$ bounds a boundary connected sum of positive-definite plumbings $\natural_{i}P_i$, and the corresponding plumbing graph is obtained from the disjoint union of the linear plumbing graphs corresponding to the summands \cite{AcetoMcCoyPark2025DefFillings}.
Note that aside from case $(1)$, each lens space in the list above does not necessarily bound a rational homology ball, thus it cannot embed individually in $\mathbb{C}P^2$. 
This could be of interest, as all results known up to now concerning embeddings of connected sums of lens spaces in $\mathbb{C}P^2$ are about lens spaces which are not obstructed from embedding individually.
\newline\indent
Further work is motivated by the following problem.

\begin{conj}[Koll\'ar \cite{Kollar}]\label{conj: kollar}
If $M$ is a Koll\'ar manifold with boundary $Y_{1} \sqcup \dots \sqcup Y_{n}$ and each $Y_i$ spherical, then:
\[
\sum\limits_{i=1}^{n} \left(1 -  \frac{1}{|\pi_{1}(Y_{i})|} \right) \leq 3.
\]
\end{conj}
Lens spaces are finite cyclic quotients of $S^3$, and hence spherical $3$-manifolds. If we consider a disjoint union of rational homology balls $B_{1} \sqcup \dots \sqcup B_{n}$ and their complement $M$ in $\mathbb{C}P^2$, we see that Koll\'ar's conjecture together with \cite[Proposition 2.4]{GollaOwensTopBMY} imply that we may only embed up to three disjoint rational homology balls bounded by lens spaces in $\mathbb{C}P^2$.
In \cite{JoParkParkMontgomeryYang}, Jo--Park--Park use lattice constraints coming from Donaldson's Diagonalisation Theorem and the $d$-invariant correction term to prove some special cases of Conjecture \ref{conj: kollar}. This motivates research of triples of rational homology balls that disjointly embed in $\mathbb{C}P^2$. In \cite{GollaOwensQHB}, Golla--Owens describe two constructions to find rational balls bounded by lens spaces that embed in a homotopy $\mathbb{C}P^2$, and that can be implemented in a computer search, producing a variety of examples.
Many triples of lens spaces, including some where at least one is of Lisca type $(2)$ or $(3)$, were found using the search implementation. The current step in understanding these embeddings is to isolate a string arising from one of these lens spaces and shed light on any possible recurrence relations among triples that could lead to a Farey-tree type of argument, as in \cite{GollaOwensQHB}. As mentioned, the constructions developed by Golla--Owens yield embeddings of triples of lens spaces into a homotopy $\mathbb{C}P^2$, therefore it would also be natural to study when this corresponds to a standard $\mathbb{C}P^2$.
\section{Tool Use}
This article and the computations presented were made entirely by the author. The software GAP \cite{GAP4} was employed for calculations as in Example \ref{ex: sporadic} and related calculations in \cite{CarforaPhD} (where it has been openly acknowledged). No AI tools were used.

    \addcontentsline{toc}{chapter}{Bibliography}
    \bibliographystyle{plainurl} 
    \bibliography{literature}

@article {Lisca2007,
    AUTHOR = {Lisca, Paolo},
     TITLE = {Lens spaces, rational balls and the ribbon conjecture},
   JOURNAL = {Geom. Topol.},
  FJOURNAL = {Geometry \& Topology},
    VOLUME = {11},
      YEAR = {2007},
     PAGES = {429--472},
      ISSN = {1465-3060,1364-0380},
   MRCLASS = {57M25 (57M12)},
  MRNUMBER = {2302495},
MRREVIEWER = {Darren\ D.\ Long},
       DOI = {10.2140/gt.2007.11.429},
       URL = {https://doi.org/10.2140/gt.2007.11.429},
}

@misc{AcetoGollaLarsonLecuona2020,
      title={Surgeries on torus knots, rational balls, and cabling}, 
      author={Paolo Aceto and Marco Golla and Kyle Larson and Ana G. Lecuona},
      year={2020},
      eprint={2008.06760},
      archivePrefix={arXiv},
      primaryClass={math.GT},
      url={https://arxiv.org/abs/2008.06760}, 
}

@article {Owens2020,
    AUTHOR = {Owens, Brendan},
     TITLE = {Smooth, nonsymplectic embeddings of rational balls in the
              complex projective plane},
   JOURNAL = {Q. J. Math.},
  FJOURNAL = {The Quarterly Journal of Mathematics},
    VOLUME = {71},
      YEAR = {2020},
    NUMBER = {3},
     PAGES = {997--1007},
      ISSN = {0033-5606,1464-3847},
   MRCLASS = {53A20 (32Q40 53D05)},
  MRNUMBER = {4142720},
MRREVIEWER = {Ahmed\ Lesfari},
       DOI = {10.1093/qmathj/haaa013},
       URL = {https://doi.org/10.1093/qmathj/haaa013},
}

@article {AcetoMcCoyPark2025DefFillings,
    AUTHOR = {Aceto, Paolo and McCoy, Duncan and Park, JungHwan},
     TITLE = {Definite fillings of lens spaces},
   JOURNAL = {J. Differential Geom.},
  FJOURNAL = {Journal of Differential Geometry},
    VOLUME = {131},
      YEAR = {2025},
    NUMBER = {1},
     PAGES = {1--41},
      ISSN = {0022-040X,1945-743X},
   MRCLASS = {57K30},
  MRNUMBER = {4947548},
       DOI = {10.4310/jdg/1755541518},
       URL = {https://doi.org/10.4310/jdg/1755541518},
}

@article {DonaldsonDiagonalization1987,
    AUTHOR = {Donaldson, Simon K.},
     TITLE = {The orientation of {Y}ang-{M}ills moduli spaces and
              {$4$}-manifold topology},
   JOURNAL = {J. Differential Geom.},
  FJOURNAL = {Journal of Differential Geometry},
    VOLUME = {26},
      YEAR = {1987},
    NUMBER = {3},
     PAGES = {397--428},
      ISSN = {0022-040X,1945-743X},
   MRCLASS = {57R19 (32C10 32G13 57N13 57R55 58E15)},
  MRNUMBER = {910015},
MRREVIEWER = {Nicholas\ Buchdahl},
       URL = {http://projecteuclid.org/euclid.jdg/1214441485},
}

@article {EvansSmith2018,
    AUTHOR = {Evans, Jonathan D. and Smith, Ivan},
     TITLE = {Markov numbers and {L}agrangian cell complexes in the complex
              projective plane},
   JOURNAL = {Geom. Topol.},
  FJOURNAL = {Geometry \& Topology},
    VOLUME = {22},
      YEAR = {2018},
    NUMBER = {2},
     PAGES = {1143--1180},
      ISSN = {1465-3060,1364-0380},
   MRCLASS = {53D35 (14J17 53D42)},
  MRNUMBER = {3748686},
MRREVIEWER = {Alexander\ Fel\cprime shtyn},
       DOI = {10.2140/gt.2018.22.1143},
       URL = {https://doi.org/10.2140/gt.2018.22.1143},
}

@article {HackingProkhorov2010,
    AUTHOR = {Hacking, Paul and Prokhorov, Yuri},
     TITLE = {Smoothable del {P}ezzo surfaces with quotient singularities},
   JOURNAL = {Compos. Math.},
  FJOURNAL = {Compositio Mathematica},
    VOLUME = {146},
      YEAR = {2010},
    NUMBER = {1},
     PAGES = {169--192},
      ISSN = {0010-437X,1570-5846},
   MRCLASS = {14J26 (14B05 14B07 14E30 14J30)},
  MRNUMBER = {2581246},
MRREVIEWER = {Mark\ Gross},
       DOI = {10.1112/S0010437X09004370},
       URL = {https://doi.org/10.1112/S0010437X09004370},
}

@article {OzsvathSzabo2003AbsolutelyGraded,
    AUTHOR = {Ozsv\'ath, Peter and Szab\'o, Zolt\'an},
     TITLE = {Absolutely graded {F}loer homologies and intersection forms
              for four-manifolds with boundary},
   JOURNAL = {Adv. Math.},
  FJOURNAL = {Advances in Mathematics},
    VOLUME = {173},
      YEAR = {2003},
    NUMBER = {2},
     PAGES = {179--261},
      ISSN = {0001-8708,1090-2082},
   MRCLASS = {57R58 (57M27)},
  MRNUMBER = {1957829},
MRREVIEWER = {Jacob\ Andrew\ Rasmussen},
       DOI = {10.1016/S0001-8708(02)00030-0},
       URL = {https://doi.org/10.1016/S0001-8708(02)00030-0},
}

@article {Hantzsche1938,
    AUTHOR = {Hantzsche, Walter},
     TITLE = {Einlagerung von {M}annigfaltigkeiten in euklidische {R}\"aume},
   JOURNAL = {Math. Z.},
  FJOURNAL = {Mathematische Zeitschrift},
    VOLUME = {43},
      YEAR = {1938},
    NUMBER = {1},
     PAGES = {38--58},
      ISSN = {0025-5874,1432-1823},
   MRCLASS = {99-04},
  MRNUMBER = {1545714},
       DOI = {10.1007/BF01181085},
       URL = {https://doi.org/10.1007/BF01181085},
}

@misc{GollaOwensQHB,
      title={The {F}arey tree and embeddings of lens spaces and rational balls in $\mathbb{CP}^2$}, 
      author={Marco Golla and Brendan Owens},
      year={2025},
      eprint={2512.09183},
      archivePrefix={arXiv},
      primaryClass={math.GT},
      url={https://arxiv.org/abs/2512.09183}, 
}

@misc{GollaOwensTopBMY,
    author = {Marco Golla and Brendan Owens },
    title = {On a class of 4-manifolds, after {Koll\'ar} -- {Unpublished}},
    journal = {},
    year = {2026}
}

@manual{GAP4,
    organization = "The GAP~Group",
    title        = "{GAP -- Groups, Algorithms, and Programming,
                    Version 4.15.1}",
    year         = 2025,
    url          = {https://www.gap-system.org},
    }

@article{HirschImbedding,
 ISSN = {0003486X, 19398980},
 URL = {http://www.jstor.org/stable/1970293},
 author = {Morris W. Hirsch},
 journal = {Annals of Mathematics},
 number = {3},
 pages = {494--497},
 publisher = {[Annals of Mathematics, Trustees of Princeton University on Behalf of the Annals of Mathematics, Mathematics Department, Princeton University]},
 title = {The Imbedding of Bounding Manifolds in Euclidean Space},
 urldate = {2026-02-11},
 volume = {74},
 year = {1961}
}

@article {EdmondsLivingston,
    AUTHOR = {Edmonds, Allan L. and Livingston, Charles},
     TITLE = {Embedding punctured lens spaces in four-manifolds},
   JOURNAL = {Comment. Math. Helv.},
  FJOURNAL = {Commentarii Mathematici Helvetici},
    VOLUME = {71},
      YEAR = {1996},
    NUMBER = {2},
     PAGES = {169--191},
      ISSN = {0010-2571,1420-8946},
   MRCLASS = {57R40 (57N13 57N35)},
  MRNUMBER = {1396671},
MRREVIEWER = {Ian\ Hambleton},
       DOI = {10.1007/BF02566415},
       URL = {https://doi.org/10.1007/BF02566415},
}

@article{LiscaParma22,
author = {Lisca, Paolo and Parma, Andrea},
title = {On {S}tein rational balls smoothly but not symplectically embedded in $\mathbb {CP}^2$},
journal = {Bulletin of the London Mathematical Society},
volume = {54},
number = {3},
pages = {949-960},
doi = {https://doi.org/10.1112/blms.12607},
url = {https://londmathsoc.onlinelibrary.wiley.com/doi/abs/10.1112/blms.12607},
year = {2022}
}

@article {LiscaSum,
    AUTHOR = {Lisca, Paolo},
     TITLE = {Sums of lens spaces bounding rational balls},
   JOURNAL = {Algebr. Geom. Topol.},
  FJOURNAL = {Algebraic \& Geometric Topology},
    VOLUME = {7},
      YEAR = {2007},
     PAGES = {2141--2164},
      ISSN = {1472-2747,1472-2739},
   MRCLASS = {57M25 (57M12 57N10)},
  MRNUMBER = {2366190},
MRREVIEWER = {Dave\ Auckly},
       DOI = {10.2140/agt.2007.7.2141},
       URL = {https://doi.org/10.2140/agt.2007.7.2141},
}

@article{Kollar,
  title={Is there a topological {B}ogomolov--{M}iyaoka--{Y}au inequality?},
  author={J{\'a}nos Koll{\'a}r},
  journal={arXiv: Algebraic Geometry},
  year={2006},
  url={https://api.semanticscholar.org/CorpusID:18236693}
}

@article {LidmanMooreVazquez,
    AUTHOR = {Lidman, Tye and Moore, Allison H. and Vazquez, Mariel},
     TITLE = {Distance one lens space fillings and band surgery on the
              trefoil knot},
   JOURNAL = {Algebr. Geom. Topol.},
  FJOURNAL = {Algebraic \& Geometric Topology},
    VOLUME = {19},
      YEAR = {2019},
    NUMBER = {5},
     PAGES = {2439--2484},
      ISSN = {1472-2747,1472-2739},
   MRCLASS = {57K10 (57K18 57R58 92E10)},
  MRNUMBER = {4023320},
MRREVIEWER = {Ioannis\ Diamantis},
       DOI = {10.2140/agt.2019.19.2439},
       URL = {https://doi.org/10.2140/agt.2019.19.2439},
}

@article {LinPin2,
    AUTHOR = {Lin, Francesco},
     TITLE = {The surgery exact triangle in {${\rm Pin}(2)$}-monopole
              {F}loer homology},
   JOURNAL = {Algebr. Geom. Topol.},
  FJOURNAL = {Algebraic \& Geometric Topology},
    VOLUME = {17},
      YEAR = {2017},
    NUMBER = {5},
     PAGES = {2915--2960},
      ISSN = {1472-2747,1472-2739},
   MRCLASS = {57M27},
  MRNUMBER = {3704248},
MRREVIEWER = {Jianfeng\ Lin},
       DOI = {10.2140/agt.2017.17.2915},
       URL = {https://doi.org/10.2140/agt.2017.17.2915},
}

@article {DonaldOwens,
    AUTHOR = {Donald, Andrew and Owens, Brendan},
     TITLE = {Concordance groups of links},
   JOURNAL = {Algebr. Geom. Topol.},
  FJOURNAL = {Algebraic \& Geometric Topology},
    VOLUME = {12},
      YEAR = {2012},
    NUMBER = {4},
     PAGES = {2069--2093},
      ISSN = {1472-2747,1472-2739},
   MRCLASS = {57M25},
  MRNUMBER = {3020201},
       DOI = {10.2140/agt.2012.12.2069},
       URL = {https://doi.org/10.2140/agt.2012.12.2069},
}

@incollection {Turaev,
    AUTHOR = {Turaev, Vladimir G.},
     TITLE = {Classification of oriented {M}ontesinos links by means of
              invariants of spin structures},
      NOTE = {Studies in topology, V},
   JOURNAL = {Zap. Nauchn. Sem. Leningrad. Otdel. Mat. Inst. Steklov.
              (LOMI)},
  FJOURNAL = {Zapiski Nauchnykh Seminarov Leningradskogo Otdeleniya
              Matematicheskogo Instituta imeni V. A. Steklova Akademii Nauk
              SSSR (LOMI)},
    VOLUME = {143},
      YEAR = {1985},
     PAGES = {130--146, 178},
      ISSN = {0373-2703},
   MRCLASS = {57M25 (57R15)},
  MRNUMBER = {806564},
MRREVIEWER = {Don\v co\ Dimovski},
}

@article {Linksignature,
    AUTHOR = {Gordon, Cameron McA. and Litherland, Richard A.},
     TITLE = {On the signature of a link},
   JOURNAL = {Invent. Math.},
  FJOURNAL = {Inventiones Mathematicae},
    VOLUME = {47},
      YEAR = {1978},
    NUMBER = {1},
     PAGES = {53--69},
      ISSN = {0020-9910,1432-1297},
   MRCLASS = {55A25},
  MRNUMBER = {500905},
MRREVIEWER = {Wilbur\ Whitten},
       DOI = {10.1007/BF01609479},
       URL = {https://doi.org/10.1007/BF01609479},
}

@article {JoParkParkMontgomeryYang,
    AUTHOR = {Jo, Woohyeok and Park, Jongil and Park, Kyungbae},
     TITLE = {Algebraic {M}ontgomery-{Y}ang problem and smooth obstructions},
   JOURNAL = {Trans. Amer. Math. Soc.},
  FJOURNAL = {Transactions of the American Mathematical Society},
    VOLUME = {378},
      YEAR = {2025},
    NUMBER = {4},
     PAGES = {2969--3003},
      ISSN = {0002-9947,1088-6850},
   MRCLASS = {14J17 (32S25 57K41)},
  MRNUMBER = {4880467},
MRREVIEWER = {Patrick\ Popescu-Pampu},
       DOI = {10.1090/tran/9364},
       URL = {https://doi.org/10.1090/tran/9364},
}

@phdthesis{CarforaPhD,
    author = {Carfora, Giulia},
    title = {Some embedding problems in 4-dimensional topology},
    school = {University of Glasgow},
    year = {2026},
    url = {https://theses.gla.ac.uk/86088/1/2026CarforaPhD.pdf}
}
\end{document}